\documentclass[11pt,reqno]{amsart} 
\usepackage{amssymb}
\usepackage{amsmath}
\usepackage{amscd}
\begin{document}
\setlength{\oddsidemargin}{0cm} \setlength{\evensidemargin}{0cm}

\theoremstyle{plain}
\newtheorem{thm}{Theorem}[section]
\newtheorem{prop}[thm]{Proposition}
\newtheorem{lemma}[thm]{Lemma}
\newtheorem{coro}[thm]{Corollary}
\newtheorem{conj}[thm]{Conjecture}
\newtheorem{cl}[thm]{Claim}

\theoremstyle{definition}
\newtheorem{Def}[thm]{Definition}
\newtheorem{ex}[thm]{Example}
\newtheorem{rem}[thm]{Remark}

\numberwithin{equation}{section} 

\title{Novikov superalgebras in low dimensions}

\author{Yifang Kang}
\address{Department of Mathematics, Tianjin University, Tianjin 300072, P.R. China}

\author{Zhiqi Chen$^*$}
\address{School of Mathematical Science \& LPMC, Nankai University,
Tianjin 300071, P.R. China}\email{chenzhiqi@nankai.edu.cn}
\thanks{$^*$ Corresponding author}

\def\shorttitle{Novikov superalgebras in low dimensions}



\begin{abstract}\noindent Novikov superalgebras are related to the quadratic
conformal superalgebras which correspond to the Hamiltonian pairs
and play fundamental role in the completely integrable systems. In
this note, we divide Novikov superalgebras into two types: $N$ and
$S$. Then we show that the Novikov superalgebras of dimension up to
3 are of type $N$.

\end{abstract}

\maketitle


\baselineskip=20pt

\section{Introduction}

A Novikov {\it super}algebra $A$ is a ${\mathbb Z}_2$-graded vector
space $A=A_0+A_1$ with a bilinear product $(u,v)\mapsto u\circ v$
for any $u\in A_i,v\in A_j$, $w\in A$ satisfying
\begin{gather} (u\circ v)\circ w-u\circ (v\circ w)=(-1)^{ij}((v\circ
u)\circ w-v\circ (u\circ w)), \\
(w\circ u)\circ v=(-1)^{ij}(w\circ v)\circ u.
\end{gather}
The even part of a given Novikov superalgebra is what is said to be
a Novikov algebra introduced in connection with the Poisson brackets
of hydrodynamic type \cite{B-N} and Hamiltonian operators in the
formal variational calculus \cite{G-D3,G-D4,G-D5,xu6,xu7}.

The supercommutator
\begin{equation} [u,v]=u\circ v-(-1)^{ij}v\circ u,\;\; \text{ for any } u\in A_i,v\in A_j
\end{equation}
makes any Novikov superalgebra $A$ a Lie superalgebra  denoted
$SLie(A)$ in what follows. The passage from a Novikov algebra $A$ to
a Lie algebra  denoted $Lie(A)$ is analogous.

The notion of Novikov superalgebra was introduced in \cite{xu1}, as
a particular case of Lie-super-admissible algebra (Gerstenhaber
called them ${\mathbb Z}_2$-graded pre-Lie algebras
\cite{gerstenhaber}). Novikov superalgebras are also related to the
quadratic conformal superalgebras \cite{xu2}.

For the notion of {\em conformal superalgebra}, see \cite{kac} (we
do not touch priority questions here). Conformal superalgebras are
related to the linear Hamiltonian operators in the
Gel$^\prime$fand-Dikii-Dorfman theory
(\cite{G-D1,G-D2,G-D3,G-D4,G-D5}) and play an important role in
quantum field theory \cite{kac} and vertex operator superalgebra
theory \cite{kac,xu3}. Recall that a conformal superalgebra
$\Re=\Re_0\oplus \Re_1$ (we denote it $(\Re,\partial,Y^+(\cdot,z))$)
is a ${\mathbb Z}_2$-graded ${\mathbb C}[\partial]$-module with a
${\mathbb Z}_2$-graded linear map $Y^+(\cdot,z):\Re\rightarrow \hom
(\Re,\Re[z^{-1}]z^{-1})$ for any $u\in \Re_i,v\in
   \Re_j$ satisfying
\begin{gather}
   Y^+(\partial u,z)=\frac{dY^+(u,z)}{dz}; \\
   Y^+(u,z)v=(-1)^{ij}Res_x \frac{e^{x\partial}Y^+(v,-x)u}{z-x}; \\
   Y^+(u,z_1)Y^+(v,z_2) \notag \\
   \qquad =(-1)^{ij}Y^+(v,z_2)Y^+(u,z_1)+Res_x
   \frac{Y^+(Y^+(u,z_1-x)v,x)}{z_2-x},
\end{gather}
where $Res_z(z^n)=\delta^n_{-1}$ for $n\in\mathbb Z$. The even part
$\Re_0$ is a {\em conformal algebra}. The definition of conformal
superalgebra in the  above generating-function form is equivalent to
the definition in \cite{kac}, where the author used a component
formula with
$Y^+(u,z)=\mathop{\sum}\limits^\infty_{n=0}\frac{u_{(n)}z^{-n-1}}{n!}$.

A {\em quadratic conformal superalgebra} is a conformal superalgebra
which is a free ${\mathbb C}[\partial]$-module over a ${\mathbb
Z}_2$-graded subspace $V$, i.e.,
\begin{equation}
\Re={\mathbb C}[\partial]V(\cong{\mathbb
C}[\partial]\otimes_{\mathbb C}V),
\end{equation}
such that the equation
\begin{equation}
Y^+(u,z)v=(w_1+\partial w_2)z^{-1}+w_3z^{-2}\; \text{ for any }
u,v\in V
\end{equation}
holds, where $w_1,w_2,w_3\in V$. A quadratic conformal superalgebra
corresponds to a Hamiltonian pair in \cite{G-D3}, which plays a
fundamental role in the theory of completely integrable systems. A
super Gel$^\prime$fand-Dorfman algebra is a ${\mathbb Z}_2$-graded
vector space $A=A_0+A_1$ with two operations $[\cdot,\cdot]$ and
$\circ$ such that $(A,[\cdot,\cdot])$ is a Lie superalgebra,
$(A,\circ)$ is a Novikov superalgebra, and the relation
\begin{equation}
[w\circ u,v]-(-1)^{ij}[w\circ v,u]+[w,u]\circ v-(-1)^{ij}[w,v]\circ
u-w\circ [u,v]=0
\end{equation}
holds for any $u\in A_i,v\in A_j,w\in A$. It was shown in \cite{xu2}
that there is a one-to-one correspondence between quadratic
conformal superalgebras and super Gel$^\prime$fand-Dorfman algebras.
It was also pointed out in \cite{xu2} that $(A,[\cdot,\cdot],\circ)$
is a super Gel$^\prime$fand-Dorfman algebra for any Novikov
superalgebra $(A,\circ)$ and the supercommutator $[\cdot,\cdot]$
relative $\circ$.

In this paper, we divide Novikov superalgebras into two types: $N$
and $S$. We show that the Novikov superalgebras of dimension $\leq
3$ are of type $N$. The full structure theory is yet to be
developed; we hope that our answer will help.

The paper is organized as follows. In \S 2, we divide Novikov
superalgebras into two types: $N$ and $S$. The Novikov superalgebras
of type $N$ are both Lie-super-admissible algebras and
Lie-admissible algebras.  Let $A=A_0+A_1$ be a Novikov superalgebra.
Then the two algebraic structures coincide if $A_1A_1=0$. We also
show that $A_1A_1\not=0$ if $A$ is of type $S$. In \S 3, we show
that the Novikov superalgebras of dimension $\leq 3$ are of type
$N$.

Throughout this paper we assume that the algebras are
finite-dimensional over $\mathbb C$. Obvious proofs are omitted.

\section{Novikov superalgebras: types $N$ and $S$}

First, we give some examples of Novikov superalgebras.
\begin{ex} Let $A$ be an associative supercommutative
superalgebra and $D$ a left $A$-module. Then $\overline{A}=D\oplus
A$ is a Novikov superalgebra if the product is defined by
\begin{equation*}(d_1+a_1)\circ (d_2+a_2)=(-1)^{ij}a_2d_1+a_1a_2\; \text{ for any }
d_1+a_1\in \overline{A}_i, d_2+a_2\in \overline{A}_j.\end{equation*}
\end{ex}

\begin{ex} Let $A$ be an associative supercommutative superalgebra.
Let $d$ be its even derivation and $c\in A_0$. Then the product
defined by
\begin{equation*}u\circ v=ud(v)+cuv\; \text{ for any }  u,v\in A\end{equation*}
determines the structure of a Novikov superalgebra on the space $A$.
\end{ex}
The former two examples are from \cite{xu2}, but we have omitted the
details inessential for the purposes of our article.

\begin{ex} Let $A=A_0+A_1$ be a ${\mathbb Z}_2$-graded vector space
with $\dim A_0=\dim A_1=1$. Let $\{e\}$ be a basis of $A_0$ and
$\{v\}$ a basis of $A_1$. Define a multiplication on $A$ by setting
$vv=e$ (we give  only non-zero products). Then $A$ is a Novikov
superalgebra and the Lie superalgebra $SLie(A)$ is defined by
$[v,v]=2e$. At the same time, $A$ is a Novikov algebra under the
above product and the Lie algebra $Lie(A)$ is abelian.
\end{ex}

\begin{ex}Let $A=A_0+A_1$ be a ${\mathbb Z}_2$-graded
vector space with $\dim A_0=1$ and $\dim A_1=2$. Let $\{e\}$ be a
basis of $A_0$ and $\{u,v\}$ be a basis of $A_1$ such that the
non-zero products are given by
\begin{equation*}
  uv=-vu=e.
\end{equation*}
Then $A$ is a Novikov superalgebra and the Lie superalgebra
$SLie(A)$ satisfies
\begin{equation*}
[x,y]=0\; \text{ for any }  x,y\in A,
\end{equation*}
which can be regarded as an abelian Lie algebra. At the same time,
$A$ is a Novikov algebra under the above product, but the Lie
algebra $Lie(A)$ is not abelian since $[u,v]=2e$.
\end{ex}

\begin{ex}
Let $A=A_0+A_1$ be a Novikov superalgebra satisfying $A_1A_1=0$.
Then the Lie superalgebra $SLie(A)$ satisfies
\begin{gather*}
   [u,v]=0,\\
   [x,y]=-[y,x],\\
   [w,[x,y]]=[[w,x],y]-[[w,y],x],\\
   [w,[u,v]]=0,
\end{gather*}
where $u,v\in A_1,w\in A$ and either $x$ or $y$ belongs to $A_0$.
Under the above product $SLie(A)$ can also be regarded as a Lie
algebra with the parity forgotten. We also have
\begin{equation*} (wu)v=0, (uv)w-u(vw)=0
\end{equation*}
for any $u,v\in A_1,w\in A$. The definitions imply that $A$ is a
Novikov algebra.
\end{ex}

\begin{Def}
Let $A=A_0+A_1$ be a Novikov superalgebra with multiplication
$(u,v)\mapsto uv$. If $A$ is also a Novikov algebra with respect to
the same product and superstructure forgotten, then $A$ is called a
Novikov superalgebra of type $N$, otherwise $A$ is said to be of
type $S$.
\end{Def}

By Definition 2.6, Novikov superalgebras in Examples 2.3--2.5 are of
type $N$. For Novikov algebras, see
\cite{B-M1,B-M2,B-M3,osborn1,osborn2,osborn3,xu4,xu5,zelmanov}.

\begin{prop} Let $A=A_0+A_1$ be a Novikov superalgebra.
If $A_1=0$ or $A_0=0$, then $A_1A_1=0$ and $A$ is of type $N$. In
particular, 1-dimensional Novikov superalgebras are of type $N$.
\end{prop}

\begin{prop}
Let $A=A_0+A_1$ be a Novikov superalgebra of type $S$. Then
$A_1A_1\not=0$.
\end{prop}

\begin{prop}
Let $A=A_0+A_1$ be a Novikov superalgebra with $\dim A_1=1$. Then
$A$ is of type $N$.
\end{prop}

\begin{rem}
Let $A=A_0+A_1$ be a Novikov superalgebra. Then $A$ is a
Lie-super-admissible algebra. If $A$ is of type $N$, then $A$ is
both a Lie-super-admissible algebra and a Lie-admissible algebra. By
Example 2.5, these two structures coincide if $A_1A_1=0$. If $A$ is
of type $N$ satisfying $A_1A_1\not=0$, then these two algebraic
structures on $A$ are different.
\end{rem}

\section{The Novikov superalgebras of $\dim\leq 3$ are of type $N$}

\begin{prop}
Any 2-dimensional Novikov superalgebra is of type $N$.
\end{prop}

\begin{prop} Let $A=A_0+A_1$ be a $3$-dimensional Novikov superalgebra of type
$S$. Then $\dim A_1=2$.
\end{prop}

\begin{proof}Assume that $A=A_0+A_1$ is a 3-dimensional Novikov superalgebra of type $S$.
By Propositions 2.7 and 2.9, we know that $\dim A_1=2$.
\end{proof}

The following discussion is  based on the modules over Novikov
algebras. The notion of a module over a Novikov algebra was
introduced in \cite{osborn1}, but a more explicit definition was
given in \cite{xu4}. A {\em module $M$ over a Novikov algebra} $A$
is a vector space endowed with two linear maps $L_M,R_M:
A\rightarrow End_{\mathbb F}(M)$ for any $x,y\in A$ satisfying
\begin{gather}
L_M(xy)=R_M(y)L_M(x), \label{1} \\
R_M(xy)-R_M(y)R_M(x)=[L_M(x),R_M(y)], \label{2}\\
R_M(x)R_M(y)=R_M(y)R_M(x), \label{3}\\
[L_M(x),L_M(y)]=L_M([x,y]). \label{4}
\end{gather}
 Let $\{e_1,e_2,\cdots,e_n\}$ be a basis of
$M$. Let $\alpha$ be a linear transformation on $M$ with the matrix
$( a_{ij})_{i,j=1}^n$ in the basis $\{e_1,e_2,\cdots,e_n\}$.

If $A$ is 1-dimensional with a basis $\{e\}$, then eqs. $(\ref{3})$
and $(\ref{4})$ are satisfied. Let $M$ be a module with a basis
$\{v_1,v_2\}$. Let $L$ and $R$ be the matrices of $L_M(e)$ and
$R_M(e)$, respectively. If $ee=0$, then by eqs. $(\ref{1})$ and
$(\ref{2})$ we have
\begin{equation}\label{ee=0}
RL=0, R^2=-LR.
\end{equation}
If $ee=e$, then by  eqs. $(\ref{1})$ and $(\ref{2})$ we have
\begin{equation}\label{ee=e}
RL=L, R^2=R+L-LR.
\end{equation}

\begin{cl}If $ee=0$, then $R^2=0$.
\end{cl}
\begin{proof} Assume that $\det R\not=0$. Then $L=0$ and $R^2=0$. It
is a contradiction. Hence $\det R=0$. If $R^2\not=0$, then $R=\left
(\begin{array}{cc}a & 0
\\ 0 & 0 \end{array}\right )$ for some $a\not=0$. It follows that $L=\left
(\begin{array}{cc}0 & 0
\\ b & c \end{array}\right )$. Then $R^2=\left
(\begin{array}{cc}a^2 & 0
\\ 0 & 0 \end{array}\right )=-LR=-\left
(\begin{array}{cc} 0 & 0
\\ ac & 0 \end{array}\right )$. Hence $a=0$. So we have $R^2=0$.
\end{proof}
\begin{cl} If $ee=e$ and $LR\not=RL$, then $L=\left (\begin{array}{cc}a & 0
\\ 0 & 0 \end{array}\right )$ for some $a\not=0$.
\end{cl}
\begin{proof}If $LR\not=RL$, then $\det L=0$. It follows that
$L=\left (\begin{array}{cc}0 & 0
\\ 1 & 0 \end{array}\right )$ or $L=\left (\begin{array}{cc}a & 0
\\ 0 & 0 \end{array}\right )$ for some $a\not=0$. For the former case, $R=\left (\begin{array}{cc}0 & 0
\\ b & 1 \end{array}\right )$. Then $R^2=R$. Hence $LR=L=RL$. It is a contradiction.
That is, $L=\left (\begin{array}{cc}a & 0
\\ 0 & 0 \end{array}\right )$ for some $a\not=0$.
\end{proof}

\begin{prop}\label{prop1}
In the above notations the classification of two-dimensional modules
over one-dimensional Novikov algebras is given in the following
table:

{\tabcolsep 0.07in \doublerulesep 0.5pt
\begin{tabular}{cccccccc} \hline \hline \hline
{\bf \footnotesize Type} & {\bf \footnotesize A} & {\bf
\footnotesize L}
    &  {\bf \footnotesize R} & {\bf \footnotesize Type} & {\bf \footnotesize A} & {\bf \footnotesize L}
    &  {\bf \footnotesize R}
  \\ \hline \hline

$T1$ & $ee=0$ & $\left(\begin{array}{cc}0 & 0
\\ 0 & 0 \end{array}\right )$ & $\left
(\begin{array}{cc}0 & 0
\\ 0 & 0 \end{array}\right )$ & $T2$ & $ee=0$  & $\left (\begin{array}{cc}1 & 0
\\ 0 & a \end{array}\right )$ & $\left
(\begin{array}{cc}0 & 0
\\ 0 & 0 \end{array}\right )$ \\

$T3$ &$ee=0$ & $\left (\begin{array}{cc}0 & 0
\\ 1 & 0 \end{array}\right )$ & $\left
(\begin{array}{cc}0 & 0
\\ 0 & 0 \end{array}\right )$ & $T4$ &
    $ee=0$ & $\left (\begin{array}{cc}1 & 0
\\ 1 & 1 \end{array}\right )$ & $\left
(\begin{array}{cc}0 & 0
\\ 0 & 0 \end{array}\right )$  \\

$T5$ &$ee=0$ & $\left (\begin{array}{cc}0 & 0
\\ a & 0 \end{array}\right )$ & $\left
(\begin{array}{cc}0 & 0
\\ 1 & 0 \end{array}\right )$ & $T6$
      & $ee=0$ & $\left (\begin{array}{cc}0 & 0
\\ a & 1 \end{array}\right )$ & $\left
(\begin{array}{cc}0 & 0
\\ 1 & 0 \end{array}\right )$ \\

$T7$ &$ee=e$ & $\left (\begin{array}{cc}0 & 0
\\ 0 & 0 \end{array}\right )$ & $\left
(\begin{array}{cc}0 & 0
\\ 0 & 0 \end{array}\right )$
     & $T8$ & $ee=e$ & $\left (\begin{array}{cc}a & 0
\\ 0 & 0 \end{array}\right )$ & $\left
(\begin{array}{cc}1 & 0
\\ 0 & 0 \end{array}\right )$ \\

$T9$ &$ee=e$ & $\left (\begin{array}{cc}a & 0
\\ 0 & a \end{array}\right )$ & $\left
(\begin{array}{cc}1 & 0
\\ 0 & 1 \end{array}\right )$ & $T10$ & $ee=e$ & $\left (\begin{array}{cc}a & 0
\\ 1 & a \end{array}\right )$ & $\left
(\begin{array}{cc}1 & 0
\\ 0 & 1 \end{array}\right )$ \\

$T11$ &$ee=e$ & $\left (\begin{array}{cc}a_1 & 0
\\ 0 & a_2 \end{array}\right )$ & $\left
(\begin{array}{cc}1 & 0
\\ 0 & 1 \end{array}\right )$ & $T12$ & $ee=e$ &  $\left (\begin{array}{cc}-1 & 0
\\ 0 & 0 \end{array}\right )$ & $\left
(\begin{array}{cc}1 & b
\\ 0 & 1 \end{array}\right )$ \\

& &$a_1\not=a_2$ & & & && $b\not=0$ \\
   \hline \hline \hline
\end{tabular}
} \\
\end{prop}
\begin{proof} Assume that $A$ is a 1-dimensional Novikov algebra and
$M$ is a 2-dimensional $A$-module  with a basis $\{v_1,v_2\}$ and
the two linear maps $L_M$ and $R_M$. Then there exists a basis
$\{e\}$ of $A$ such that $ee=0$ or $ee=e$. Let $L$ and $R$ be the
matrices of $L_M(e)$ and $R_M(e)$, respectively.

(I) $ee=0$. By Claim 3.3, $R^2=0$. Then there exists another basis
of $M$ (also denoted by $\{v_1,v_2\}$) such that $R=\left
(\begin{array}{cc}0 & 0
\\ 1 & 0 \end{array}\right )$ or $R=\left (\begin{array}{cc}0 & 0
\\ 0 & 0 \end{array}\right )$. If $R=\left (\begin{array}{cc}0 & 0
\\ 0 & 0 \end{array}\right )$,
then there exists another basis of $M$ (also denoted by
$\{v_1,v_2\}$) such that $L=\left (\begin{array}{cc}a_1 & 0
\\ 0 & a_2 \end{array}\right )$ or $L=\left (\begin{array}{cc}a & 0
\\ 1 & a \end{array}\right )$. If $L=\left
(\begin{array}{cc}a & 0
\\ 1 & a \end{array}\right )$, replacing $e$ by $\frac{e}{a}$ and $v_2$ by
$\frac{v_2}{a}$ for $a\not=0$, we can take $a=1$. If $R=\left
(\begin{array}{cc}0 & 0
\\ 1 & 0 \end{array}\right )$, then by the first part of
eq.~(\ref{ee=0}), we have $L=\left (\begin{array}{cc}0 & 0
\\ a & b \end{array}\right )$. Replacing $e$ by $\frac{e}{b}$ and $v_2$ by
$\frac{av_2}{b}$ for $b\not=0$, we can set $b=1$.

(II) $ee=e$. If $RL\not=LR$, then $L=\left (\begin{array}{cc}a & 0
\\ 0 & 0 \end{array}\right )$ for some $a\not=0$ by Claim 3.4. Furthermore
$R=\left (\begin{array}{cc}1 & b
\\ 0 & c \end{array}\right )$ for some $b\not=0$. By the second part of eq. (\ref{ee=e}),
we have $a=-c=-1$. If $LR=RL$, then $R^2=R$ by eq. (\ref{ee=e}).
Hence there exists another basis of $M$ (also denoted by
$\{v_1,v_2\}$) such that $R=\left (\begin{array}{cc}0 & 0
\\ 0 & 0 \end{array}\right )$, or $R=\left (\begin{array}{cc}1 & 0
\\ 0 & 0 \end{array}\right )$, or $R=\left (\begin{array}{cc}1 & 0
\\ 0 & 1 \end{array}\right )$. By the first part of
eq.~(\ref{ee=e}), if $R=\left (\begin{array}{cc}0 & 0
\\ 0 & 0 \end{array}\right )$, then $L=RL=\left (\begin{array}{cc}0 & 0
\\ 0 & 0 \end{array}\right )$. Similarly, if $R=\left (\begin{array}{cc}1 & 0
\\ 0 & 0 \end{array}\right )$, then $L=\left (\begin{array}{cc}a & b
\\ 0 & 0 \end{array}\right )$. So $b=0$ by $LR=RL$. If $R=\left (\begin{array}{cc}1 & 0
\\ 0 & 1 \end{array}\right )$, then there exists
another basis of $M$ (also denoted by $\{v_1,v_2\}$) such that
$R=\left (\begin{array}{cc}1 & 0
\\ 0 & 1 \end{array}\right )$ and $L=\left
(\begin{array}{cc}a_1 & 0
\\ 0 & a_2 \end{array}\right )$, where $a_1\not=a_2$, or $L=\left (\begin{array}{cc} a & 0
\\ 0 & a \end{array}\right )$, or $L=\left (\begin{array}{cc}a & 0
\\ 1 & a \end{array}\right )$.\end{proof}

\begin{prop} Let $A=A_0+A_1$ be a 3-dimensional Novikov superalgebra satisfying $A_1A_1\not=0$.
If $\dim A_1=2$, then $A_0A_0=A_0A_1=A_1A_0=0.$
\end{prop}
\begin{proof} By assumption,
$A_1$ is a 2-dimensional $A_0$-module. For the type $T2$ in
Proposition~\ref{prop1},
\begin{gather*}
(v_1v_1)v_1=0\Longrightarrow v_1v_1=0,\\
(v_2v_1)v_1=0\Longrightarrow v_2v_1=0,\\
a\not=0,(v_1v_2)v_2=0\Longrightarrow v_1v_2=0.\\
a=0,(ev_1)v_2=-(ev_2)v_1\Longrightarrow v_1v_2=0,\\
a\not=0,(v_2v_2)v_2=0\Longrightarrow v_2v_2=0,\\
a=0,(v_2v_2)v_1=-(v_2v_1)v_2=0\Longrightarrow v_2v_2=0.
\end{gather*}
That is, $A_1A_1=0$. We similarly show that $A_1A_1=0$ for the other
cases, except $T1$. Since $A_1A_1\not=0$, we see that
$A_0A_0=A_0A_1=A_1A_0=0.$
\end{proof}

\begin{prop}
Any 3-dimensional Novikov superalgebra is of type $N$.
\end{prop}
\begin{proof}Assume that $A=A_0+A_1$ is a 3-dimensional Novikov superalgebra of type $S$.
By Propositions 2.8 and 3.2, $A_1A_1\not=0$ and $\dim A_1=2$. Then
$A_0A_0=A_0A_1=A_1A_0=0$ by Proposition 3.6. It follows that $A$ is
of type $N$, which is a contradiction. \end{proof}

\begin{thm} The Novikov superalgebras of dimension up to 3 are of type $N$.
\end{thm}
\begin{proof} The theorem follows from Propositions 2.7, 3.1 and 3.7.
\end{proof}

\begin{rem} Since the Novikov superalgebras of type $N$ are essentially Novikov
algebras, there is a method to give the classification of the
Novikov superalgebras of type $N$ based on the classification of
Novikov algebras. That is, we look for a grading for any Novikov
algebra. But then we need to do it case-by-case. Moreover, the
classification in higher dimensions is also an open problem, in
particular, the classification of 4-dimensional Novikov algebras has
not been finished yet. Then this method does not work for the
Novikov algebras of dimension $>3$.
\end{rem}

\section*{Acknowledgments}
We are grateful of the referees and the editors for their valuable
comments and suggestions. Also we would like to acknowledge
Professor Bai C M for his helpful suggestions.

{}

\label{lastpage}


\begin{thebibliography}{99}
\small
\bibitem{B-M1}
\textsc{Bai C M} and \textsc{Meng D J}, The Classification of
Novikov Algebras in Low Dimensions, {\it J.Phys.A:Math.Gen.} {\bf
34}(2001), 1581--1594.

\bibitem{B-M2}
\textsc{Bai C M} and \textsc{Meng D J}, On the Realization of
Transitive Novikov Algebras, {\it J.Phys.A: Math.Gen.} {\bf
34}(2001), 3363--3372.

\bibitem{B-M3}
\textsc{Bai C M} and \textsc{Meng D J}, The Realization of
Non-transitive Novikov Algebras, {\it J.Phys.A: Math.Gen.} {\bf
34}(2001), 6435--6442.

\bibitem{B-N}
\textsc{Balinskii A A} and \textsc{Novikov S P}, Poisson Brackets of
Hydrodynamic Type, Frobenius Algebras and Lie Algebras, {\it Soviet
Math.Dokl.} {\bf 32}(1985), 228--231.

\bibitem{G-D1}
\textsc{Gel$^\prime$fand I M} and \textsc{Diki L A}, Asymptotic
Behaviour of the Resolvent of Strum-Liouville Equations and the
Algebra of the Korteweg-deVries Equations, {\it Russian Math.
Surveys.} {\bf 30}(1975), 77--113.

\bibitem{G-D2}
\textsc{Gel$^\prime$fand I M} and \textsc{Diki L A}, A Lie Algebra
Structure in a Formal Variational Calculation, {\it
Funct.Anal.Appl.} {\bf 10}(1976), 16--22.

\bibitem{G-D3}
\textsc{Gel$^\prime$fand I M} and \textsc{Dorfman I Y}, Hamiltonian
Operators and Algebraic Structures Related to Them, {\it
Funktsional. Anal. i Prilo\v zen.} {\bf 13}(1979), 13--30.

\bibitem{G-D4}
\textsc{Gel$^\prime$fand I M} and \textsc{Dorfman I Y}, Schouten
Brackets and Hamiltonian Operators, {\it Funktsional. Anal. i
Prilo\v zen.} {\bf 14}(1980), 71--74.

\bibitem{G-D5}
\textsc{Gel$^\prime$fand I M} and \textsc{Dorfman I Y}, Hamiltonian
Operators and Indinite Dimensional Lie Algebras, {\it Funktsional.
Anal. i Prilo\v zen.} {\bf 14}(1981), 23--40.

\bibitem{gerstenhaber}
\textsc{Gerstenhaber M}, The Cohomology Structure of an Associative
Ring, {\it Ann.Math.} {\bf 78}(1963), 267--288.

\bibitem{kac}
\textsc{Kac V G}, Vertex Algebras for Begginers, University Lectures
Series, {\bf 10}, American Math. Soc.Porvidence, 1996.

\bibitem{osborn1}
\textsc{Osborn J M}, Novikov Algebras, {\it Nava J.Algebra Geom.}
{\bf 1}(1992), 1--14.

\bibitem{osborn2}
\textsc{Osborn J M}, Simple Novikov Algebras with an Idempotent,
{\it Commun.Algebra.} {\bf 20}(1992), 2729--2753.

\bibitem{osborn3}
\textsc{Osborn J M}, Infinite Dimensional Novikov Algebras of
Characteristic 0, {\it J.Algebra.} {\bf 167}(1994), 146--167.

\bibitem{xu1}
\textsc{Xu X P}, Variational Calculus of Supervariables and Related
Algebraic Structures, {\it J.Algebra}, {\bf 223}(2000), 396--437.

\bibitem{xu2}
\textsc{Xu X P}, Quadratic Conformal Superalgebras, {\it J.Algebra},
{\bf 231}(2000), 1--38.

\bibitem{xu3}
\textsc{Xu X P}, Introduction to Vertex Operators Superalgebras and
Their Modules, Klvwer Academic, Dordercht/Boston/London, 1998.

\bibitem{xu4}
\textsc{Xu X P}, On Simple Novikov Algebras and Their Irreducible
Modules, {\it J.Algebra.} {\bf 185}(1996), 905--934.

\bibitem{xu5}
\textsc{Xu X P}, Novikov-Poisson Algebras, {\it J.Algebra.} {\bf
190}(1997), 253--279.

\bibitem{xu6}
\textsc{Xu X P}, Hamiltonian Operators and Associative Algebras with
a Derivation, {\it Lett.Math.Phys.} {\bf 33}(1995), 1--6.

\bibitem{xu7}
\textsc{Xu X P}, Hamiltonian Superoperators, {\it J.Phys
A:Math.Gen.} {\bf 28}(1995), 1681--1698.

\bibitem{zelmanov}
\textsc{Zel$^\prime$manov E I}, On a Class of Local Translation
Invariant Lie Algebras, {\it Soviet Math.Dokl.} {\bf 35}(1987),
216--218.
\end{thebibliography}
\end{document}